\documentclass[center 10pt]{article}
\usepackage{mathbbold}
\usepackage{comment}
\usepackage{mathrsfs}
\usepackage{amsmath}
\usepackage{amsfonts}
\usepackage{amssymb}
\usepackage{color}
\usepackage{array}
\usepackage{graphicx}
\usepackage{amsmath}
\usepackage{amssymb}
\usepackage{fancyhdr}
\usepackage{mathrsfs}
\usepackage{mathbbold}
\usepackage{amsfonts}
\usepackage{wasysym}
\usepackage{amscd}
\usepackage{bbm}
\usepackage{color}
\usepackage{array}
\usepackage{graphicx}
\usepackage{graphics}
\newtheorem{Theorem}{Theorem}[section]

\newtheorem{Remark}[Theorem]{Remark}

\def\Int2{\int_{-1}^1}
\def\Intt3{\int_{-{\pi\over2}}^{\pi\over2}}
\usepackage{makecell}
\title{Moment-based cosh-Hilbert Inversion and Its Applications in Single-photon Emission Computed Tomography}
\author{Shousheng Luo, Jiansheng Yang,
 Tie Zhou\\ LMAM, School of Mathematical Sciences, Peking University,
Beijing 100871, China}
\date{}
\begin{document}
\large
\maketitle
{\bf Abstract: }
%We proposed a moment-based image reconstruction method in single-photon emission computed tomography (SPECT).
The inversion of cosh-Hilbert transform (CHT) is one of the most crucial steps for single-photon emission computed tomography with uniform attenuation from truncated projection data. Although the uniqueness of the CHT inversion had been proved \cite{Noo2007}, there is no exact and analytic inverse formula so far. Several approximated inversion algorithms of the CHT had been developed \cite{Noo2007}\cite{You2007}. In this paper, we proposed a new numerical moment-based inversion algorithm.
%\begin{center}
%\Large{求解cosh-Hilbert变换的矩方法及其在单光子发射计算机断层成像中的应用}\\
%\normalsize
%\vspace{0.3cm}
%罗守胜\hspace{1cm}杨建生 \hspace{1cm}  周铁\\
%\vspace{0.3cm}
%中国, 北京, 北京大学数学科学学院，数学及其应用教育部重点实验室, 100871
%%Moment-based cosh-Hilbert Inversion and Its Applications in Single-photon Emission Computed Tomography
%\end{center}
%{\bf 摘要： }
%%We proposed a moment-based image reconstruction method in single-photon emission computed tomography (SPECT).
%当衰减系数是常数时，cosh-Hilbert变换（CHT）的逆运算是截断数据单光子发射计算机断层成像的关键步骤之一。尽管CHT变换
%的唯一性已经被证明\cite{Noo2007}，但是目前并没有精确的解析求解公式。文献\cite{Noo2007}\cite{You2007}的作者已提出了
%不同的求解CHT的近似算法。本文提出了一种新的基于矩的CHT求解算法。
%%The inversion of cosh-Hilbert transform (CHT) is one of the most crucial step for single-photon emission computed tomography with uniform attenuation from truncated projection data. Although the uniqueness of the CHT inversion had been proved \cite{Noo2007}, there is no exact and analytic inverse formula so far. Several approximated inversion algorithms of the CHT had been developed \cite{Noo2007}\cite{You2007}. In this paper, we proposed a new numerical moment-based inversion algorithm.
%\\
%{\bf 关键词：}cosh-Hilbert变换；矩方法；SPECT；截断数据
\section{Introduction}
Single-photon emission computed tomography(SPECT)
 is a non-invasive diagnostic technology in nuclear medicine.
It is used to image the physiological function of various organs with the help of radiopharmaceutical, a biochemical
molecular labeled by radioisotope.
Radiopharmaceutical is induced into a patient, and
 the emitted gamma photons, attenuated by the body, are measured by a gamma camera rotating around it.
The task of SPECT reconstruction is to
estimate the radiopharmaceutical distribution from  the measured data.
\\
\indent From the mathematical point of view, the two-dimensional(2D) SPECT reconstruction problem is to inverse
the attenuated Radon transform(aRt) of $p$ \cite{Natterer2001},
\begin{eqnarray}
g(s,\phi)=R_\mu{p}(s,\phi)=\int_\mathbb{R}p(s\theta+t\theta^\perp)e^{-\int_t^\infty \mu(s\theta+\tau\theta^\perp){d\tau}}dt,\label{aRT}
\end{eqnarray}
where $p(x)$ and $\mu(x)$ denote the distribution of radiopharmaceutical and the attenuation coefficient of gamma photon at $x=(x_1,x_2)$. Further, $\theta = (\cos\phi, \sin\phi)$ and $\theta^\perp = (-\sin\phi, \cos\phi)$ denote two perpendicular vectors.
Although $p$ is a function with a compact support $\Omega$ in practice,
%Although the
%integrand in (\ref{aRT}) is zero outside a bounded interval, with respect to $s$,
the integral (\ref{aRT}) is written over $(-\infty, +\infty)$ for convenience.
\\
\indent
 If $\mu=0$, the aRt formula (\ref{aRT}) is the well-known Radon transform(Rt), which is the mathematical basis of computed tomography(CT). We can inverse Rt by the filtering back-projection(FBP) formula \cite{Natterer2001} from 180$^\circ$ projection data due to the relation that $(R{p})(-s,\phi+\pi) = (R{p})(s,\phi)$, where $R$ denotes the Radon transform.
If $\mu\neq 0$,
the filtering back-projection formula for aRt was firstly discovered by the author of \cite{Novikov2002}
 (Latter, a simpler deviation was illustrated in \cite{Natterer_a2001}).
% and Kunyansky implemented the formula numerically \cite{Kunyansky2001}.
\\
\indent
In this paper, we are interested in a special case of aRt, which assumes that the attenuation map $\mu$ is an uniform constant in the convex domain including the support of $p$. With the assumption that $\mu(x)=\mu_0$(constant) inside $\Omega$ and $\mu(x)=0$ outside $\Omega$, the reconstruction of $p$
is equivalent to solving the exponential Radon transform(eRt) of $p$ (See \cite{Metz1980}\cite{Natterer2001} for details)
\begin{eqnarray}
g(s,\phi)=R_{\mu_0}p(s,\phi)=\int_\mathbb{R}p(s\theta+t\theta^\perp)e^{\mu_0t}dt\label{eRt}
\end{eqnarray}
%For SPECT, the eRt is applicable to brain imaging \cite{You2009}.
%\\
\indent
Since the exponential Radon transform does not have the parity property of the Radon transform, i.e.
$(R_{\mu_0}{p})(-s,\phi+\pi) \neq (R_{\mu_0}{p})(s,\phi)$ for $\mu_0\neq 0$, it is a long-standing opinion that $2\pi$ data is necessary for the reconstruction of $p$.
However, accurate reconstruction
from non-truncated projections known only over 180$^\circ$ was shown to be possible for any value
of $\mu_0$ \cite{Noo2001,Pan2002}.
%the authors of \cite{Noo2001}\cite{Noo2007}\cite{}
%proved that $p$ could be reconstructed uniquely with $180^\circ$ non-truncated projections.
%, and proposed an iterative reconstruction method \cite{Noo2001}.
%The results  were extended to various $\pi$-scheme data acquisitions,
%which allow the $180^\circ$ of data to be distributed over a union of disjoint intervals.
\\
\indent
Because of
practical constraints due to the imaging hardware, scanning geometry, or ionizing radiation exposure,
 the projection data are not theoretically
sufficient for exact image reconstruction. %Insufficient data problems occur quite frequently.
%The insufficient data problem can take many forms,
One of the important case is the truncated data that the projection data  $g(\phi,s)$ are only known for a limited range of $s$.
For example, if $\Omega$ is the centered disc of radius $R$, the
projection data is truncated as long as $g(\phi, s)$ is not known for all $s\in[-R,R]$.
\\
\indent
An explicit formula was suggested in \cite{Rull2004}, in which the author
reduced the reconstruction of 2D SPECT to solving a one-dimensional integral equation
called cosh-weighted Hilbert transform(CHT) by the weight-differential backprojection  of eRt.
This method further was developed further in \cite{Noo2007}, in which the author
proved a similar result of the classic Radon transform. Recently, this result was extended to developed
the reconstruction of three-dimensional(3D) SPECT with uniform attenuation\cite{You2009}. More importantly,
this method can be used to conduct the reconstruction of region-of-interesting from truncated projection data
\cite{Noo2007,Rull2004}.  However, there is no explicit formula for CHT so far.
Several numerical algorithms of the CHT have been proposed \cite{Noo2007,You2007,Rull2004}.
%For the classic Radon transform, we can reconstruct $p$ partially
%from truncated projection data by using backprojected filtration
%formula \cite{Noo_ct2004}.
%For the exponential Radon transform, it had been pointed out that the WDB of eRt could reduce
%the inverse of eRt to solve the CHT, which could be used to conduct the reconstruction of region-of-interesting from truncated projection data
%\cite{Noo2007,Rull2004}.
% However, there is no explicit formula for CHT so far.
%In this paper, we propose another numerical algorithm to inverse the CHT.
%\\
%\indent For the SPECT reconstruction from truncated data,
%the method proposed in \cite{Rull2004} change inversion of eRt
%into the inversion of one-dimensional integraton transform called
%cosh-weighted Hilbert transform (CHT) of $p$. This method provide
%the possibility to reconstruct part of $p$ from $180^\circ$ truncated
%projection data. The authors of \cite{Noo2007} developed this method further,
%and the author gave a approximate inversion formula. \\
\\
\indent
In this paper, we propose a moment-based method to solve the CHT numerically. In this proposed approach,
the solution of CHT is represented as the result of Tricomi formula for Hilbert transform followed by a correction term, which is related to the moments of the underlying function.
 The numerical simulations show that the proposed method
can solve the CHT efficiently for a large range of the attenuation constant $\mu_0$.
%We first propose to inverse cosh-weighted Hilbert transform by conjugate gradient method.
%In addition, a approximate with regularization technique is proposed.
\\
\indent The rest of this paper is organized as follow. We review the derivation of cosh-weighted Hilbert transform, and present the moment-based method  in section \ref{section_Num}.
Numerical experiments are shown in section \ref{num_res}. Some discussions and conclusions are given
in section \ref{con_dis}. And the appendix A  is devoted to some numerical techniques.
%\end{comment}
%\end{itemize}
%%%%%%%%%%%%%%%%%%%%%%%%
\section{Cosh-weighted Hilbert Transform and the Proposed Method}\label{section_Num}
%\begin{enumerate}
%\item power method;
%\item conjugate method.
%\end{enumerate}
This section consists of two parts.
We firstly review the derivation of cosh-weighted Hilbert transform(CHT) from the exponential Radon transform\cite{Noo2007}\cite{Rull2004}.
Then, the moment-based method is given.
%\\
%\indent
%{\red Equation (13) reduces the reconstruction problem to solving a 1D integral equation with
%a convolution kernel. However, note that $p$ is not readily found from $b$ using a simple 1D
%deconvolution, because the convolution kernel $\frac{\cosh(\mu_0x_2)}{x_2\pi }$ grows exponentially with
%increasing $y$.}
\subsection{Review of Cosh-weighted Hilbert Transform}
The basic idea is similar to differential back-projection of Radon transform ($\mu_0=0$), but with a weight function.
%%%%%%%%%%%%%%%%%
%\begin{equation}
%\int_0^\pi{e^{\mu_0x\cdot{\theta}^\perp}}(E_{\mu_0})(\phi,x\cdot\theta)d\theta\notag
%\end{equation}
%%%%%%%%%%%%
% Eq 4
Let $b(x)$ be defined as
\begin{equation}
b({x})=\int_0^\pi{e^{-\mu_0 x\cdot{\theta}^\perp}}(R_{\mu_0}^\prime{p})(x\cdot\theta,\phi)d\phi\label{Dbp},
\end{equation}
where  the derivative is respect to the first variable of $R_{\mu_0}p$, i.e.,
$(R_{\mu_0}^\prime{p})(s,\phi)={{\partial}\over{\partial{s}}}(R_{\mu_0}{p})(s,\phi)$, and $x\cdot\theta$ denotes
the standard inner product of $\mathbb{R}^2$.
%Replacing $(R_\mu_0{p})(s,\phi)$ by its integral expression (\ref{eq2}),
By the definition of eRt (\ref{eRt}),
interchanging the partial
derivative with the integration and applying the chain rule, we  can obtain
\begin{equation}
(R_{\mu_0}^\prime{p})(s,\phi)=\int_\mathbb{R}{e^{\mu_0 t^\prime}}\theta\cdot\nabla{p}(s\theta+t^\prime\theta^\perp)dt^\prime\notag.
\end{equation}
Here, we assume that $p$ is continuously differentiable. Since
$
x = (x\cdot\theta)\theta + (x\cdot\theta^\perp)\theta^\perp
$,
we have that
\begin{eqnarray}
(R_{\mu_0}^\prime{p})({x}\cdot\theta,\phi)&=&\int_\mathbb{R}{e^{\mu_0{t}^\prime}}\theta\cdot\nabla{p}((x\cdot\theta)\theta+t^\prime\theta^\perp)dt^\prime\notag\\
                          &=&\int_\mathbb{R}{e^{\mu_0{t}^\prime}}\theta\cdot\nabla{p}({x}+(t^\prime-{x}\cdot\theta^\perp)\theta^\perp)dt^\prime\notag\\
&=&e^{\mu_0{x}\cdot\theta}\int_\mathbb{R}{e^{\mu_0{t}}}\theta\cdot\nabla{p}(x+t\theta^\perp)dt \notag                        \\
&=&-e^{\mu_0{x}\cdot\theta}\int_\mathbb{R}{{e^{\mu_0{t}}}\over{t}}{{\partial}\over{\partial\phi}}p(x+t\theta^\perp)dt,\label{dbp1}
 \end{eqnarray}
where we use a variable replacement of $t=t^\prime-{x}\cdot\theta$ for the third equality, and the last equality is implied by
\begin{eqnarray}
{{\partial}\over{\partial\phi}}p(x+t\theta^\perp)
&=&\nabla{p(x+t\theta^\perp)}\cdot(t{{\partial}\over{\partial\phi}}\theta^\perp)\notag\\
%\theta^\perp=(-\sin\phi,\cos\phi)
&=&-t\theta\cdot\nabla{p(x+t\theta^\perp)}\notag.
\end{eqnarray}
%\end{slide}
%%%%%%%%%%%%%%%%%%%%%%%%%
%\begin{slide}
%%%%%%%%%%
% Eq. 9
Inserting (\ref{dbp1}) into (\ref{Dbp}) and interchanging the order of intergation, we have that
\begin{eqnarray}
b({x})&=&-\int_\mathbb{R}{{e^{\mu_0{t}}}\over{t}}p(x+t\theta^\perp)\left|_0^\pi{dt}\right.\notag\\
    &=&-\int_\mathbb{R}{{e^{\mu_0{t}}}\over{t}}\left[p(x_1,x_2-t)-p(x_1,x_2+t)\right]dt\notag\\
    &=&-\int_\mathbb{R}{{e^{\mu_0{t}}}\over{t}}p(x_1,x_2-t)dt+\int_\mathbb{R}{{e^{\mu_0{t}}}\over{t}}p(x_1,x_2+t)dt\label{eq_repac}\\
    &=&-2\pi\int_\mathbb{R}{{\cosh{\mu_0(x_2-x_2^\prime)}}\over{\pi(x_1-x_2^\prime)}}p(x_1,x_2^\prime)dx_2^\prime,\label{eq3}
\end{eqnarray}
where we use variable replacements of $x_2^\prime=x_2-t$ and $x_2^\prime=x_2+t$ for the first and second terms of (\ref{eq_repac}), respectively. The above integrals (\ref{eq_repac}) and (\ref{eq3}) are understood in the sense of Cauchy principal value. In fact, we are only interested in the reconstruction of $p$ at $(x_1, x_2) \in \Omega$,
a bounded and convex region outside which $p$ is known to be zero. Denoting by $(x_1,L(x_1))$ and $(x_1,U(x_1))$ the end points of the intersection of
$\Omega$ with the line parallel to the $x_2$-axis through the point $(x_1,x_2) \in\Omega$, we have
\begin{eqnarray}
b(x_1,x_2)=-2\pi\int_{L(x_1)}^{U(x_1)}{{\cosh{\mu_0(x_2-x_2^\prime)}}\over{\pi(x_2-x_2^\prime)}}p(x_1,x_2^\prime)dx_2^\prime
\end{eqnarray}
If $\Omega$ is a centered disc of radius $R$, we have $U(x_1)=-L(x_1)=\sqrt{R^2-x_1^2}$.
\\
\indent %Therefore, the inversion of eRt is equivalent to solve the CHT along each vertical line.
%In addition to the CHT function, we also know the values $R_{\mu_0}(x_1,0)$ and $R_{\mu_0}(-x_1,\pi)$.
Therefore, the SPECT reconstruction problem is changed into solving the so-called finite cosh-weighted Hilbert transform along each vertical line with $(R_{\mu_0}p)(x_1,0)$ and $(R_{\mu_0}p)(-x_1,\pi)$ known.
\\
%%%%%%%%%%%%%%%%%%
\indent
%In the following two subsections, we will discuss the inversion o
For simplicity,
 we focus on the following standard
form: reconstruct a function $f(t)$ with its support $\text{supp}{(f)}\subset(-1,1)$
 from  its finite cosh-weighted Hilbert transform
\begin{eqnarray}
h_{\mu_1}(\tau)   &=&\int_{-1}^{1}{{\cosh{\mu_1(\tau-t)}}\over{\pi(\tau-t)}}f(t)dt,\text{\indent\ } |\tau|<1
\label{eq_cosh}
%\mu_0    &=&\mu_0{{U_x-L_x}\over{2}}\\
%f(\tau)&=&p(x,{{U_x-L_x}\over{2}}\tau+{{U_x+L_x}\over{2}})\\
%h(t)   &=&{{1}\over{\pi(L_x-U_x)}}b(x,{{U_x-L_x}\over{2}}t+{{U_x+L_x}\over{2}})\\
%m_\mu_0  &=&{{2}\over{(L_x-U_x)}}(e^{-\mu_0{{U_x+L_x}\over{2}}}E_{\mu_0}f(0,x)+e^{\mu_0{{U_x+L_x}\over{2}}}E_{\mu_0}f(\pi,-x))
\end{eqnarray}
with $c_{\mu_1}=\int_{-1}^1f(t)\cosh({\mu_1t})dt$ known. We can normalize the general form
of (\ref{eq3}) as the standard form by using the following affine transformation
%\begin{eqnarray}
%\mu_1&=&\mu_0\frac{U(x_1)-L(x_1)}{2},\\
%f(\tau)& = &p(x_1,{U(x_1)+L(x_1)\over2}+\tau{U(x_1) - L(x_1)\over2}),\\
%h(t)& = &\frac{1}{\pi(L(x_1) - U(x_1))}b(x_1,\frac{U(x_1) + L(x_1)}{2}+\tau\frac{U(x_1) - L(x_1)}{2}).
%\end{eqnarray}
\begin{eqnarray}
\mu_1&=&\hat{d}\mu_0,\notag\\
f(\tau)& = &p(x_1,\hat{c}+\hat{d}\tau),\notag\\
h_{\mu_1}(t)& = &-{1\over2\hat{d}\pi}b(x_1,\hat{c}+\hat{d}\tau),\notag
\end{eqnarray}
where $\hat{c}={U(x_1)+L(x_1)\over2}$ and $\hat{d}=\frac{U(x_1)-L(x_1)}{2}$.
Further, we have
\begin{eqnarray}
c_{\mu_1}=\frac{1}{\hat{d}}\left[e^{-\hat{c}\mu_0}(R_{\mu_0}p)(x_1,0)
+e^{\hat{c}\mu_0}(R_{\mu_0}p)(-x_1,\pi)\right].\notag
\end{eqnarray}
If $\mu_1=0$, by the Tricomi formula \cite{Tricomi}, we have
\begin{eqnarray}
f(t)\sqrt{1-t^2}&=&{{c_0}\over{\pi}}-\int_{-1}^1{\sqrt{1-s^2}\over{\pi(t-s)}}h_0(s)ds\label{tricomi1}
%c_{0}             &=&\int_{-1}^1f(s)ds\notag.
\end{eqnarray}
Then, we can obtain $f$ for $|t|<1$ by dividing $\sqrt{1-t^2}$ on both sides of (\ref{tricomi1}).
\subsection{Moment-based Method}
In this subsection, we focus on the inverse of the finite cosh-weighted Hilbert transform (\ref{eq_cosh}). For the
sake of reference, we introduce some notations firstly.
%For simplicity, we introduce the following notes: $K_\mu_0=({H_\mu_0,C_\mu_0})$
\begin{eqnarray}
h_m(t)&=&{1\over\pi}\int_{-1}^1(t-\tau)^{m-1}f(\tau)d\tau,\notag\\
T_n(t)&=&{1\over\pi}\Int2{\sqrt{1-s^2}\over{t-s}}s^nds,\notag\\ T_i^j&=&{1\over\pi}\Int2\frac{t^j}{\sqrt{1-t^2}}T_i(t)dt,\notag
\end{eqnarray}
and the $m$-th moment $c_m=\int_{-1}^1f(t)t^mdt$ of $f$. Obviously,
we have that $T_n(-t)=(-1)^{n+1}T_n(t)$. \\
%Firstly, we have that that $T_n(-t)=(-1)^{n+1}T_n(t)$.
 %we have the following theorems \ref{Theorem1} and \ref{Theorem2}.
%In order to maintain the clearly frame of this work, the proof of them are present in the appendix A.
%%%%%%%%%%%%%%%%%%%%%%%%%% Theorem 1%%%%%%%%%%%%
%By the Taylor expansion of $\cosh$ function and the definition of CHT,
\indent
By the Taylor expansion of the function $\cosh\mu_1t$, we have
\begin{eqnarray}
h_{\mu_1}(t)&=&{1\over\pi}\int_{-1}^1\left[{1\over{t-\tau}}f(\tau)+\sum\limits_{k=1}^\infty{\mu_1^{2k}\over{(2k)!}}(t-\tau)^{2k-1}f(\tau)\right]d\tau\notag\\
        &=&h_0(t)+\sum\limits_{k=1}^\infty{\mu_1^{2k}\over{(2k)!}}{h_{2k}}(t),\label{taylor1}
\end{eqnarray}
 and
\begin{eqnarray}
c_{\mu_1}   &=&\int_{-1}^1\cosh(\mu_1{t})f(t)dt\notag\\
        %&=&\sum_{k=0}^\infty\frac{\mu_1^{2k}}{(2k)!}\int_{-1}^1t^{2k}f(t)dt\notag\\
        &=&c_0+\sum\limits_{k=1}^\infty{\mu_1^{2k}\over{(2k)!}}{c_{2k}}\label{taylor2}.
\end{eqnarray}
%where we define $0!=1$.
Therefore, $h_0(t)$ and $c_0$ can be represented by
%\begin{eqnarray}
$
h_0(t)=h_{\mu_1}(t)-\sum\limits_{k=1}^\infty{\mu_1^{2k}\over{(2k)!}}{h_{2k}}(t)$, $
c_0   =c_{\mu_1}-\sum\limits_{k=1}^\infty{\mu_1^{2k}\over{(2k)!}}{c_{2k}}.
%\end{eqnarray}
$
Applying the Tricomi formula for the finite Hilbert transform\cite{Tricomi},
we have
\begin{eqnarray}
&&f(t){\sqrt{1-t^2}}\notag\\
%&=&{{c_0}\over{\pi}}-{1\over\pi}\Int2{{\sqrt{1-\tau^2}}\over{t-\tau}}h_0(\tau)d\tau\notag\\
    &=&{1\over\pi}(c_{\mu_1}-\sum\limits_{k=1}^\infty{\mu_1^{2k}\over{(2k)!}}{c_{2k}})
    -{{1}\over{\pi}}\Int2{{\sqrt{1-\tau^2}}\over{t-\tau}}\left(h_{\mu_1}(\tau)-\sum\limits_{k=1}^\infty{\mu_1^{2k}\over{(2k)!}}{h_{2k}}(\tau)\right)d\tau\notag\\
    %&=&f_0(t)
    %-{{1}\over{\pi}}\sum\limits_{k=1}^\infty{\mu_1^{2k}\over{(2k)!}}{c_{2k}}+{1\over{\pi^2}}\Int2{{\sqrt{1-\tau^2}}\over{t-\tau}}\Int2\sum\limits_{k=1}^\infty{\mu_1^{2k}\over{(2k)!}}(\tau-s)^{2k-1}f(s)dsd\tau\notag\\
    &=&f_{\mu_1}(t)
    -{{1}\over{\pi}}\sum\limits_{k=1}^\infty{\mu_1^{2k}\over{(2k)!}}{c_{2k}}+{1\over{\pi^2}}\sum\limits_{k=1}^\infty{\mu_1^{2k}\over{(2k)!}}\Int2{{\sqrt{1-\tau^2}}\over{t-\tau}}\Int2(\tau-s)^{2k-1}f(s)dsd\tau\notag\\
    &=&f_{\mu_1}(t)
    -{{1}\over{\pi}}\sum\limits_{k=1}^\infty{\mu_1^{2k}\over{(2k)!}}{c_{2k}}+{1\over{\pi^2}}\sum\limits_{k=1}^\infty{\mu_1^{2k}\over{(2k)!}}\Int2f(s)\Int2{{\sqrt{1-\tau^2}}\over{t-\tau}}(\tau-s)^{2k-1}d\tau ds,\notag\\
&=&f_{\mu_1}(t)
    -{{1}\over{\pi}}\sum\limits_{k=1}^\infty{\mu_1^{2k}\over{(2k)!}}{c_{2k}}+{1\over{\pi^2}}\sum\limits_{k=1}^\infty{\mu_1^{2k}\over{(2k)!}}\sum\limits_{l=0}^{2k-1}
    A_{2k-1}^l\Int2f(s)(-s)^{l}\Int2{{\sqrt{1-\tau^2}}\over{t-\tau}}\tau^{2k-1-l}d\tau ds\notag\\
    &=&f_{\mu_1}(t)
    -{{1}\over{\pi}}\sum\limits_{k=1}^\infty{\mu_1^{2k}\over{(2k)!}}{c_{2k}}+{1\over{\pi}}\sum\limits_{k=1}^\infty{\mu_1^{2k}\over{(2k)!}}\sum\limits_{l=0}^{2k-1}A_{2k-1}^l
    (-1)^lc_{l}T_{2k-1-l}(t),\label{eq_represnent}
\end{eqnarray}
\large
where
$%\begin{eqnarray}
f_{\mu_1}(t)={1\over\pi}(c_{\mu_1}-\Int2{{\sqrt{1-\tau^2}}\over{t-\tau}}h_{\mu_1}(\tau)d\tau)
$ and
 $A_n^m=\frac{n!}{m!(n-m)!}$. The equation (\ref{eq_represnent}) provides us a relationship between the moments of $f$ and $f_{\mu_1}$, which tells us that we can reconstruct $f$ accurately if the moments of
 $f$ are known.
\\
\indent Define $d_n=\int_{-1}^{1}{t^n\over\sqrt{1-t^2}}f_{\mu_1}(t)dt$ firstly.
Multiplying $\frac{t^{2i}}{\sqrt{1-t^2}}$( $i=0,1,2,\cdots, \infty$) on both sides of (\ref{eq_represnent}),
and integrating on the interval $(-1,1)$, we have that
\small
\begin{eqnarray}
c_{2i}&=&\Int2f(t)t^{2i}dt\notag\\
      &=&\Int2\frac{f_{\mu_1}(t)t^{2i}}{\sqrt{1-t^2}}dt-\left(\Int2\frac{t^{2i}}{\sqrt{1-t^2}}dt\right)\sum\limits_{k=1}^\infty{\mu_1^{2k}\over{(2k)!}}c_{2k}+
      {1\over{\pi}}\sum\limits_{k=1}^\infty{\mu_1^{2k}\over{(2k)!}}\sum\limits_{l=0}^{2k-1}A_{2k-1}^{l}
(-1)^lc_{l}\left({1\over\pi}\Int2{T_{2k-1-l}(t)t^{2i}\over{\sqrt{1-t^2}}}dt\right)\notag\\
      &=&d_{2i}-B_{2i}\sum\limits_{k=1}^\infty{\mu_1^{2k}\over{(2k)!}}c_{2k}+
      \sum\limits_{k=1}^\infty{\mu_1^{2k}\over{(2k)!}}\sum\limits_{l=0}^{k-1}A_{2k-1}^{2l}
    T_{2(k-l)-1}^{2i}c_{2l},\label{eq_even1}
        %&=&d_{2i}-B_{2i}\sum\limits_{k=1}^\infty{\mu_1^{2k}\over{(2k)!}}{c_{2k}}+
    %\sum_{k=0}^\infty[\sum_{l=k}^\infty{\mu_1^{2(l+1)}\over{(2(l+1))!}}{A}_{2l+1}^{2k}T_{2(l-k)+1}^{2i}]c_{2k}\label{eq_even1}
\end{eqnarray}
\large
where $B_{2i}={1\over\pi}\Int2\frac{t^{2i}}{\sqrt{1-t^2}}dt={1\over\pi}
\int_{-\frac{\pi}{2}}^{\frac{\pi}{2}}\sin^{2i}\alpha{d\alpha}=\frac{(2i-1)!!}{(2i)!!}$ for $i>0$, and $B_0=1$. The second equality is implied by the uniform convergence of the series, and the third equality holds because $T_{2k}(t)$ is an odd function for all $k$ by the definition of $T_n(t)$.
\\
\indent Similarly, multiplying $\frac{t^{2i-1}}{\sqrt{1-t^2}}$( $i=1,2,\cdots, \infty$) on both sides of (\ref{eq_represnent}), and integrating on the interval $(-1,1)$, we have that
\begin{eqnarray}
c_{2i-1}&=&d_{2i-1}-\sum\limits_{k=1}^\infty{\mu_1^{2k}\over{(2k)!}}\sum_{l=1}^{k}
      A_{2k-1}^{2l-1}T_{2(k-l)}^{2i-1}c_{2l-1}\label{eq_odd1}.
%    \sum_{k=1}^\infty[\sum_{l=k}^\infty{\mu_1^{2l}\over{(2l)!}}{A}_{2l-1}^{2k-1}T_{2(l-k)}^{2i}]c_{2k-1}.\label{eq_odd1}
\end{eqnarray}
Here, we use the fact that  $T_{2k-1}(t)$ is an even function for all $k$.
Reformulating the equations (\ref{eq_even1}) and (\ref{eq_odd1}),
we can derive two infinite linear systems with respect to the odd and even
moments of $f$, respectively.
\begin{eqnarray}
Qc_{even}&=&d_{even},\label{eq_even}\\
Pc_{odd}&=&d_{odd},\label{eq_odd}
%c_{2i-1}&=&\sum\limits_jP_{ij}c_{2j-1},i=1,2,\cdots\\
%c_{2i}&=&\sum\limits_jQ_{ij}c_{2j},i=0,1,\cdots,
\end{eqnarray}
%where for all $i,j\geq 1$, %$i=1,2,\cdots,+\infty$ and
where $c_{even}=(c_0,c_2,\cdots,c_{2k},\cdots)$, $c_{odd}=(c_1,c_3,\cdots,c_{2k+1},\cdots)$, and
$d_{even}$ and $d_{odd}$ have the same meanings. Further,
\begin{eqnarray}
Q_{ij}&=&\delta_{ij}+{\mu_1^{2j}\over{(2j)!}}{{(2i-1)!!}\over{(2i)!!}}-\sum\limits_{k=j}{\mu_1^{2k}\over{(2(k+1))!}}A_{2k+1}^{2j}T_{2k+1-2j}^{2i},\\
P_{ij}&=&\delta_{ij}+\sum\limits_{k=j}^\infty{\mu_1^{2k}\over{(2k)!}}A_{2k-1}^{2j-1}T_{2(k-j)}^{2i-1}.
\end{eqnarray}
%We have the following linear systems with respect to the odd and even moments of $f$, respectively.
%
%\begin{eqnarray}
%%\sum\limits_jP_{ij}c_{2j-1}&=&d_{2i-1},i=1,2,\cdots\label{sys_1}\\
%Pc_{odd}&=d_{odd},\label{sys_1}\\
%%\sum\limits_jQ_{ij}c_{2j}&=&d_{2i},i=0,1,\cdots\label{sys_2},
%Qc_{even}&=d_{even},\label{sys_2}
%\end{eqnarray}
%%and $Q$ and $P$ are invertible operators.
%where Let $c_{odd}=(c_1,c_3,\cdots)$, $c_{even}=(c_0,c_2,\cdots)$, and $d_{odd}=(d_1,d_3,\cdots),d_{even}=(d_0,d_2,\cdots)$ with
%$d_i=\Int2\frac{t^i}{\sqrt{1-t^2}}f_{\mu_1}(t)dt$, $P$ and $Q$ are two infinite operators which elements are defined as
%\begin{eqnarray}
%P_{ij}&=&\delta_{ij}+\sum\limits_{k=j}^\infty{\mu_1^{2k}\over{(2k)!}}A_{2k-1}^{2j-1}T_{(2k-j)}^{2i-1},\label{eq_even}\\
%Q_{ij}&=&\delta_{ij}+{\mu_1^{2j}\over{(2j)!}}{{(2i-1)!!}\over{(2i)!!}}-\sum\limits_{k=j}^{\infty}{\mu_1^{2k}\over{(2(k+1))!}}A_{2k+1}^{2j}T_{2k+1-2j}^{2i}.\label{eq_odd}
%\end{eqnarray}
Here $\delta_{ij}$ denotes the Kronecker delta.
Due to the uniqueness of CHT \cite{Noo2007}, the equation systems (\ref{eq_even})(\ref{eq_odd}) can be solved uniquely.
%\\
%\indent
% In addition to the reconstruction formula (\ref{synth_1}), we have two other reconstruction methods theoretically
% if the moments of $f$ are known.
\begin{Remark}
Although we can approximate $f$ by polynomials if the moments of $f$ are known, the simulations show
that the reconstructed images are dominated by distortions because of the computed errors and the high frequency of $f$.
\end{Remark}
%The equation (\ref{eq_represnent}) provides us a relationship between the moments of $f$ and $f_{\mu_1}$.
Therefore,
we can synthesize $f(t)\sqrt{1-t^2}$ by formula (\ref{eq_represnent}), then divide $\sqrt{1-t^2}$ to get $f$ finally.
However, in fact we are only able to get a finite number of moments of $f$
by truncating the equation (\ref{eq_even}) and (\ref{eq_odd}).
Fortunately, we only need a few low moments in computation, because of the rapid decrease of $|c_j|$ with the increase of $j$.
\\
\indent
 If we truncate the first $2M+1$ terms of (\ref{taylor1}) and (\ref{taylor2}), we can obtain a finite approximation of $f(t)\sqrt{1-t^2}$
\begin{eqnarray}
f(t){\sqrt{1-t^2}}
    &\approx&f_{\mu_1}(t)
    -{{1}\over{\pi}}\sum\limits_{k=1}^M{\mu_1^{2k}\over{(2k)!}}{c_{2k}}+{1\over{ \pi}}\sum\limits_{k=1}^{2M}{\mu_1^{2k}\over{(2k)!}}\sum\limits_{l=0}^{2k-1}
    (-1)^lA_{2k-1}^lc_{l}T_{2k-1-l}(t).\label{finit_synth}
\end{eqnarray}
Further, we have finite approximations
of the operator equations (\ref{eq_even}) and (\ref{eq_odd}).
\begin{eqnarray}
\hat{P}\hat{c}_{odd}\approx\hat{d}_{odd}\label{fsys_odd}\\
\hat{Q}\hat{c}_{even}\approx\hat{d}_{even}\label{fsys_even}
\end{eqnarray}
where
%$\hat{c}_{even},\hat{c}_{odd},\hat{d}_{even},\hat{d}_{odd}$ are the
%truncation of $\hat{c}_{even},\hat{c}_{odd},\hat{d}_{even},\hat{d}_{odd}$, respectively.
$\hat{d}_{odd}=(d_1,d_3,\cdots,d_{2M-1})$ for $1\leq i,j\leq M$
%with $k=\lfloor\frac{M+1}{2}\rfloor$
and $\hat{d}_{even}=(d_0,d_2,\cdots,d_{2M})$. %with $k=\lfloor{M\over2}\rfloor$.
Further, $\hat{c}_{even}$ and
$\hat{c}_{odd}$ have the same meaning as $\hat{d}_{even}$ and $\hat{d}_{odd}$, respectively.
Further,
\begin{eqnarray}
\hat{P}_{ij}&=&\delta_{ij}+\sum\limits_{k=j}^M{\mu_1^{2k}\over{(2k)!}}A_{2k-1}^{2j-1}T_{2(k-j)}^{2i-1},\\
\hat{Q}_{ij}&=&\delta_{ij}+{\mu_1^{2j}\over{(2j)!}}{{(2i-1)!!}\over{(2i)!!}}-
\sum\limits_{k=j}^M{\mu_1^{2k}\over{(2(k+1))!}}A_{2k+1}^{2j}T_{2k+1-2j}^{2i}.
\end{eqnarray}
 If
the integer of $M$ is large enough,
we have that the solutions to the equations (\ref{fsys_odd}) and (\ref{fsys_even}) are unique by the
continuity of the spectral of operators, and
we can obtain the first $2M+1$ moments of
$f$ by solving the finite linear equation systems (\ref{fsys_odd}) and (\ref{fsys_even}).
%Here $\lfloor{r}\rfloor$
%denotes the largest integer that is not larger than $r$.
\\
\indent
Based on the discussions above, we can summarize the reconstruction algorithm of SPECT as the following steps: \begin{enumerate}
\item  performing the weight-differential backprojection (\ref{Dbp}) to compute the CHT along each vertical line intersected with the
       image domain, which is used to reducing the reconstruction problem into solving one-dimensional integral equations(CHT);
\item computing $f_{\mu_1}$ by Tricomi formula, $d_i=\Int2\frac{t^i}{\sqrt{1-t^2}}f_{\mu_1}(t)dt$ for $i=0,1,2,\cdots,2M$ further, and the system matrices $\hat{P}$ and $\hat{Q} $, then solving the linear systems (\ref{eq_even}) and (\ref{eq_odd}) to get the moments $c_i$ of $f$ for $i=0,2\cdots,2M$;
\item synthesizing $f(t)\sqrt{1-t^2}$ by (\ref{finit_synth}), then dividing $\sqrt{1-t^2}$ to get $f(t)$ for $|t|<1$.
\end{enumerate}
\section{Numerical Results}\label{num_res}
\begin{figure}[h]
\centering
\begin{tabular}{cc}
\includegraphics[width=5cm]{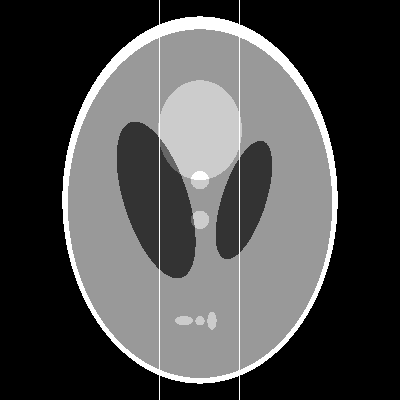}
&
\includegraphics[width=5cm]{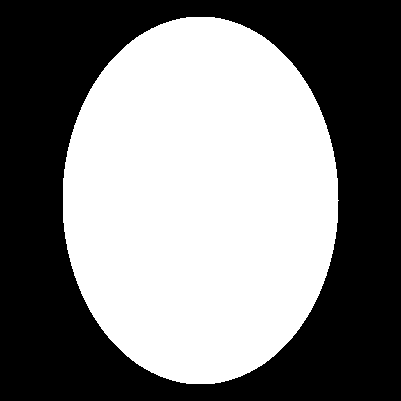}\\
(a)&(b)
\end{tabular}
\caption{Left: activity distribution (a). Right: attenuation map domain (b).}\label{fig-phan}
\end{figure}

%Figure~\ref{white} is an example of Pop-Art.
%\begin{figure}[!hbp]
%\makebox[\textwidth]{\framebox[5cm]{\rule{0pt}{5cm}}}
%\caption{Five by Five in Centimetres.} \label{white}
%\end{figure}

\begin{figure}
\centering
\begin{tabular}{cccc}
\includegraphics[width=4cm]{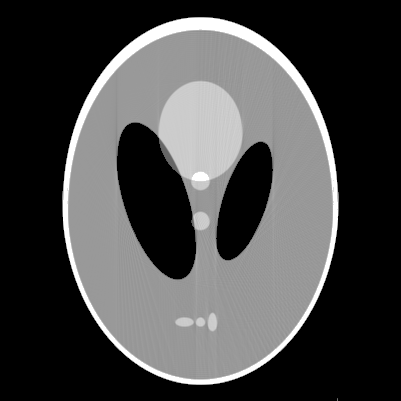}
\includegraphics[width=4cm]{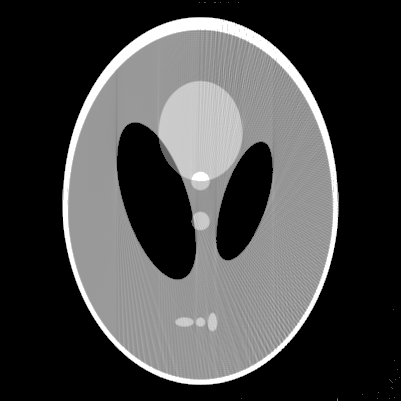}\\
%&\includegraphics[width=4cm]{cosh-data/clear/400-1000/u4}\\
%\includegraphics[width=4cm]{cosh-data/clear/trun401/u0}
\includegraphics[width=4cm]{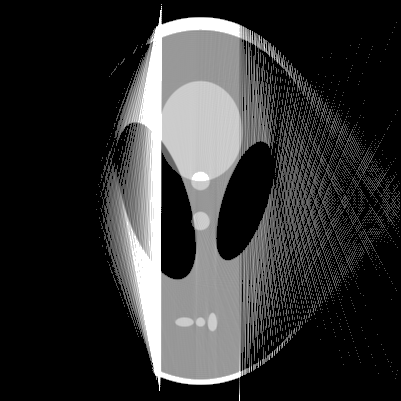}
\includegraphics[width=4cm]{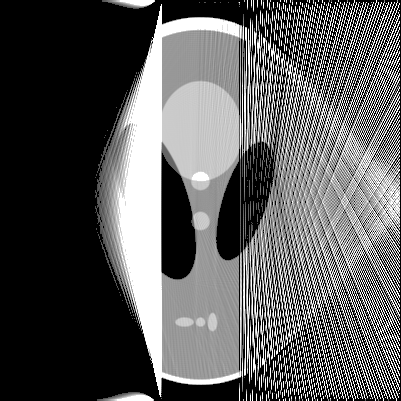}
%&\includegraphics[width=4cm]{cosh-data/clear/trun401/u4}
\end{tabular}
\caption{Images reconstructed from the noise-free data by the proposed algorithm.
}\label{fig-perfect}
\end{figure}

\begin{figure}[h]
\centering
\begin{tabular}{ccc}
\includegraphics[width=4cm]{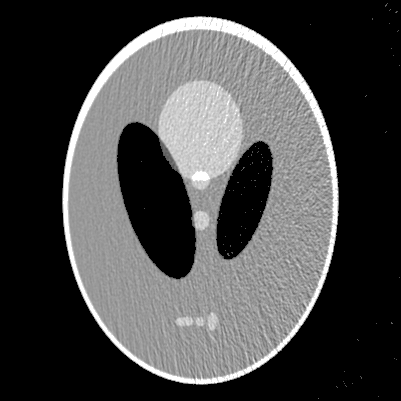}
\includegraphics[width=4cm]{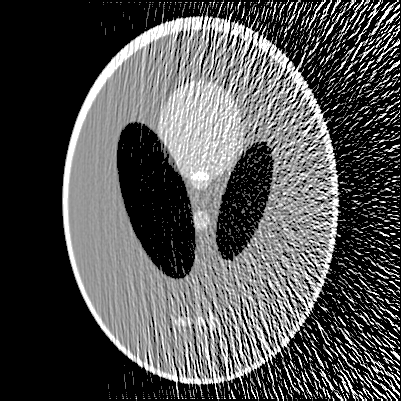}\\
%&\includegraphics[width=4cm]{cosh-data/noise/cosh(401)/u4}\\
%\includegraphics[width=4cm]{cosh-data/noise/trun(401)/9/u0}
\includegraphics[width=4cm]{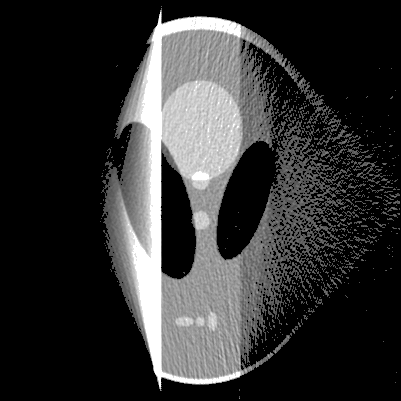}
\includegraphics[width=4cm]{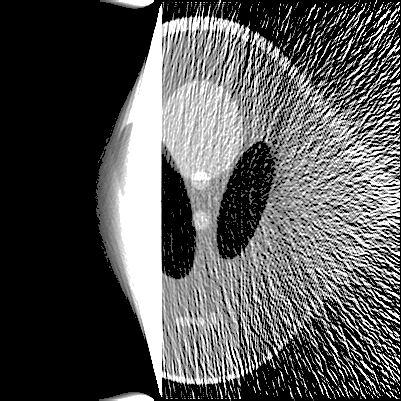}
%&\includegraphics[width=4cm]{cosh-data/noise/cosh(401)/u4}
\end{tabular}
\caption{Images reconstructed from the noisy projection data by the proposed algorithm.
%not passing through the rectangular box in figure \ref{fig-phan} was discarded.
%The images in top row are reconstructed from noise-free projection data, and these in the bottom row from
%Poisson noised projection data corresponding to a total count of $10^9$ photons.
%From left to right, in each row, the attenuation coefficient $\mu_1$ is successively 1.5 cm$^{-1}$ and 3 cm$^{-1}$ and 4 cm$^{-1}$.
}\label{fig-noise}
\end{figure}

\begin{table}[h]
\centering
\caption{ Definition of the ellipses forming the SPECT version of the Shepp-Logan phantom.}\label{tab-phan}
\begin{tabular}{ccccc}
\hline
\Xhline{0.5pt}
Centre(cm)    &1st axis (cm) &2nd axis (cm) &Polar angle($^\circ$) & Intensity\\
\hline
(0, 0)          &0.69           &0.92           &0   &0.5\\
(0, -0.0184)    &0.6624         &0.874          &0   &-0.2\\
(0.22, 0)       &0.31           &0.11           &72  &-0.2\\
(-0.22, 0)      &0.41           &0.16           &108 &-0.2\\
(0, 0.35)       &0.21           &0.25           &0   &0.1\\
(0, 0.1)        &0.046          &0.046          &0   &0.1\\
(0, -0.1)       &0.046          &0.046          &0   &0.1\\
(-0.08, -0.605) &0.046          &0.023          &0   &0.1\\
(0, -0.605)     &0.023          &0.023          &0   &0.1\\
(0.06, -0.605)  &0.203          &0.046          &0   &0.1\\
\Xhline{1pt}
\end{tabular}
\end{table}
%\begin{table}
%\centering
%\caption{ Definition of the ellipses forming the SPECT version of the Shepp-Logan phantom.}\label{tab-phan}
%\begin{tabular}{ccccc}
%\hline
%\Xhline{0.5pt}
%Centre(cm)    &1st axis (cm) &2nd axis (cm) &Polar angle($^\circ$) & Intensity\\
%\hline
%(0, 0)        &6.9           &9.2           &0   &0.5\\
%(0, -0.184)   &6.624         &8.74          &0   &0.2\\
%(2.2, 0)      &3.1           &1.1           &72  &0.2\\
%(-2.2, 0)     &4.1           &1.6           &108 &0.2\\
%(0, 3.5)      &2.1           &2.5           &0   &0.1\\
%(0, 1)        &0.46          &0.46          &0   &0.1\\
%(0, -1)       &0.46          &0.46          &0   &0.1\\
%(-0.8, -6.05) &0.46          &0.23          &0   &0.1\\
%(0, -6.05)    &0.23          &0.23          &0   &0.1\\
%(0.6, -6.05)  &0.23          &0.46          &0   &0.1\\
%\Xhline{1pt}
%\end{tabular}
%\end{table}
In this section, we will validate the performance of the proposed method by numerical experiments.
 For this purpose, a SPECT version of the Shepp-Logan phantom(see figure \ref{fig-phan}(a), \cite{Noo2007}) was used under various conditions.
 %(see \cite{Noo2007} for the description of the ellipses in details).
Table \ref{tab-phan} gives a description of the ellipses that form this phantom.
% The support of attenuation map is defined by figure \ref{fig-phan}(b), and in experiments
 The attenuation map we used is equal to a constant inside the ellipse shown in figure \ref{fig-phan}(b). Although we used different values in the ellipse for different experiments to test the stability of the proposed methods, we only illustrate the reconstructed images in figures \ref{fig-perfect} and \ref{fig-noise} for $\mu_1=1.5$ and $3$, which are sufficiently large
for the medical applications of SPECT(see \cite{Noo2007} for details).
\\
\indent
%Perfect projection data were created via exponential Radon transform formula
%(\ref{eRt}). % , which included tissue attenuation, but neither scatter nor blurring.
In the experiments, perfect projection data at 1000 views sampled over $[0,\pi]$ evenly with 400 rays
are created via the eRt formula (\ref{eRt}).
In order to test the stability and robust of the proposed method about noise, the data were Poisson
noised by the procedure used in \cite{Noo2007}.
The reconstructions were discretized in a $ 400\times 400$ grid.
%in order to simulate quantum noise in the data, the following procedure
%was implemented.
%\begin{itemize}
%\item The projections were scaled (multiplied by a constant factor) so that the total count was a predefined integer,
%% easured?value in the data set was equal to 50,
%\item each point value in the data set was then replaced by a random realization of a Poisson
%variate with a mean equal to that value.
%\end{itemize}
%Each projection contains {\red 129} rays.

\indent
The proposed method was performed for various computer-simulated data, including perfect, noisy,
non-truncated and truncated data.
Figure \ref{fig-perfect} displays the reconstructed images from perfect data,
while figure \ref{fig-noise} displays the reconstruction from noised data for $\mu_1=1.5$ and $3$. The images in top row are reconstructed from non-truncated projection data, while these in the bottom row are reconstructed from the
truncated data such that any measurement corresponding to a line
not passing through the rectangular box in figure \ref{fig-phan}(a) was discarded.

%From left to right, in each row, the attenuation coefficient $\mu_1$ is successively 1.5 cm$^{-1}$ and 3 cm$^{-1}$.
\begin{figure}
\centering
\includegraphics[width=14cm]{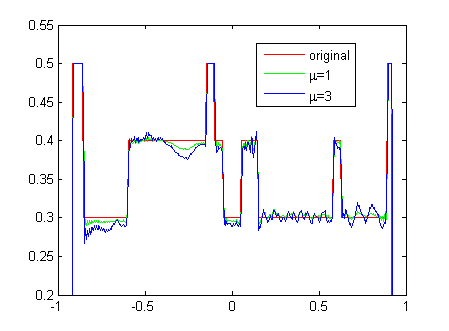}
\caption{Profiles of the images reconstructed from perfect non-truncated projection data in figure \ref{fig-perfect}.}\label{fig-prof}
\end{figure}

\indent
 From figures \ref{fig-perfect} and \ref{fig-noise}, we can see that the presented algorithms can carry out the reconstruction of SPECT efficiently.
 As our expectation, the distortions in the reconstructions raise with the increase of the attenuation coefficient $\mu_1$. The profiles of the reconstructions in figure \ref{fig-perfect}
at $x_1=0$ are displayed in figure \ref{fig-prof}.
Similar to \cite{Noo2007}, we have the same
observation that the noise in the reconstructions appears strongly spatial variant, for which the reason is the weight variant in (\ref{Dbp})(see \cite{Noo2007} for details).

%\indent In our simulations, the projection angles are over $[0,\pi]$,
%which rotates the detectors on the left sides {\red figure 2}. In the
%acquisition procedure the distance between any pixel on the right part
%of the phantom and the detectors are larger than that between pixels on the left part and the detectors.
%\input{tab}
%%%%%%%%%%%%%%%%%%%%%%%%%%%%%%%%%
\section{Conclusions and discussions}\label{con_dis}
In this work,
a numerical method for the CHT is presented. The numerical experiments show that the proposed method can
conduct the SPECT reconstruction for practical attenuation constant approximately. On the other hand, we
can observe that there are great distortions with the increasing of attenuation constant.
Though numerical methods can inverse the CHT approximately, great efforts are needed to seek a analytic inversion formula and more robust numerical methods for larger attenuation constant $\mu_1$.
\section*{Acknowledgements}
This work is supported by the
National Basic Research Program of China (2011CB809105) and NSF grants of China (61121002, 10990013).
The authors are grateful for helpful discussions with Professor Haomin Zhou(School of Mathematics,
Georgia Institute of Technology).
%We also give great thanks to Professor Jiansheng Yang
%(School of Mathematical Sciences, Peking University) for the useful suggestions.
\section*{Appendix}
In this appendix, we give the analytic formulae for the computation of $T_i(t)$ and $T_i^j$.
And we will discuss the numerical method of $d_i$, which is very important for the
accurate reconstruction. %$T_i(t)$ and $T_i^j$ have the following recursion relations.
\\
{\bf Computation of $T_n(t)$: }
Firstly, we have \cite{Tricomi}
\begin{eqnarray}
%T_{-1}(t)&=&1\\
T_0(t)&=&{1\over\pi}\Int2\frac{\sqrt{1-\tau^2}}{t-\tau}d\tau=t.
\end{eqnarray}
By this result, we can compute $T_n(t)$($n>0$) by recursive relationship
\begin{eqnarray}
%T_0(s)&=&{1\over\pi}\int_{-1}^1{{\sqrt{1-t^2}}\over{s-t}}dt=s\notag\\
T_n(s)&=&{1\over\pi}\int_{-1}^1{t^n{\sqrt{1-t^2}}\over{s-t}}dt\notag\\
      &=&{1\over\pi}\int_{-1}^1t^{n-1}{{s-(s-t)}\over{s-t}}\sqrt{1-t^2}dt\notag\\
      &=&{1\over\pi}\int_{-1}^1t^{n-1}{{1}\over{s-t}}\sqrt{1-t^2}dt-{1\over\pi}\int_{-1}^1t^{n-1}\sqrt{1-t^2}dt\notag\\
      &=&sT_{n-1}(s)-S_{n-1}\label{eq_recu},
\end{eqnarray}
where
$
S_{n}=\left\{\begin{array}{ll}
0&n \textrm{  is odd}\\
B_{2n}-B_{2n+2}& n \textrm{  is even}
\end{array}
\right.,
$
and
$B_{2n}={1\over\pi}\int_{-{\pi\over2}}^{\pi\over2}\sin^{2n}(\alpha)d\alpha=\frac{(2n-2)!!}{(2n)!!}$.
%Because $T_0(t)=t$ is an odd function, we have that the parity of $T_n(t)$ is contrary to odevity of $n$ by
%the recursive formula (\ref{eq_recu}), i.e. $T_n(-t)=(-1)^{n+1}T_n(t)$.
\\
{\bf Computation of $T_i^j(t)$: }
 Firstly, we have that $T_n(-t)=(-1)^nT(t)$ by the definition of $T_n(t)$. Therefore,
we have $T_{2n}^{2k}=0,T_{2n+1}^{2k+1}=0$ for $n,k\geq0$, and $T_{0}^{2k+1}={{(2k+1)!!}\over{(2k+2)!!}}=B_{2k+2}$ because $T_0(t)=t$.
By using the recursive formula (\ref{eq_recu}) of $T_n(t)$, we have
\begin{eqnarray}
T_{2n}^{2k+1}&=&{1\over\pi}\Int2{{s^{2k+1}\over\sqrt{1-s^2}}T_{2n}(s)}ds\notag\\%={1\over\pi}\Int2{{s^{2k+1}\over\sqrt{1-s^2}}(s^2T_{2n-2}(s)-sS_{2n-2})}ds\notag\\
%     &=&T_{2n-2}^{2k+3}-S_{2n-2}*B_{2k+2}\\
     &=&{1\over\pi}\Int2{{s^{2k+1}\over\sqrt{1-s^2}}(sT_{2n-1}(s))}ds\notag\\
     &=&T_{2n-1}^{2k+2},
\end{eqnarray}
and
\begin{eqnarray}
T_{2n+1}^{2k}&=&{1\over\pi}\Int2{s^{2k}\over\sqrt{1-s^2}}T_{2n+1}(s)ds\notag\\
             &=&{1\over\pi}\Int2{s^{2k}\over\sqrt{1-s^2}}(sT_{2n}(s)-S_{2n})ds\notag\\
             %&=&{1\over\pi}\Int2{s^{2k+1}\over\sqrt{1-s^2}}T_{2n}(s)ds-S_{2n}{1\over\pi}\Int2{s^{2k}\over\sqrt{1-s^2}}ds\notag\\
             &=&T_{2n}^{2k+1}-S_{2n}B_{2k}.
\end{eqnarray}
%{\bf The generated steps for $T_{odd}$ and $T_{even}$}\\
%It is very important to compute $d_i$ precisely in equation {\red()} for for the inversion (\ref{synth_1}) to compute $c_i$ precisely.
{\bf Computation of $d_i$: }
In order to get the unknown function $f(t)$ for $|t|<1$ finally, we should divide $\sqrt{1-t^2}$ on both sides of (\ref{finit_synth}),
which implies that
tiny errors near the two ends $-1$ and 1 must result in large deviations from $f(t)$.
Therefore, we have to compute $c_i$ and $d_i$ precisely to suppress the errors on two ends.
Because the computation formula for each $d_i$ can be seen as the $\frac{1}{\sqrt{1-t^2}}$ weighted integral of $f_{\mu_1}(t)t^i$, we use the Gauss-Chebyshev quadrature formula\cite{quadrature} in this paper.
\begin{eqnarray}
{1\over\pi}\int_{-1}^{1}\frac{s(t)}{\sqrt{1-t^2}}dt\approx{\pi\over n}\sum_{i=1}^ns(q_i), \label{eq_comp_d}
\end{eqnarray}
 where the Gauss-Chebyshev nodes of (\ref{eq_comp_d}) are given by
 $q_i=\cos{\frac{2i-1}{2n}\pi}$. %  are the Chebyshev's point, i.e $\cos(t_i)=0$.}
%In the simulations of this work, the values at $q_i$ by spline interpolation scheme from
%uniformly discrete points on the interval $(-1,1)$.
\bibliographystyle{ieeetr}
\bibliography{cosh_Hilbert6}
\end{document}